\documentclass[12pt]{article} 
\usepackage{amsmath}
\usepackage{amssymb}
\usepackage{cite}
\newsymbol  \blackbox 1004
\newcommand{\eh}{\hfill}\newlength{\sperr}

\newenvironment{proof}{{\settowidth{\sperr}{\bf\rm
Proof}%
\par\addvspace{0.3cm}\noindent\parbox[t]{1.3\sperr}
{\bf\rm P\eh r\eh o\eh o\eh f\eh }%
}}{\nopagebreak\mbox{}
$\blackbox$\par\addvspace{0.3cm}}

\def\nn{\nonumber}

\def\diag{{\mathrm{diag}}}
\def\a{\alpha}
\def\b{\beta}
\def\g{\gamma}
\def\Ga{\Gamma}

\def\vk{\varkappa}
\def\Lam{\Lambda}
\def\s{\sigma}
\def\la{\lambda  }
\def\om{\omega}
\def\Om{\Omega}
\def\S{\Sigma}
\def\t{\theta}
\def\Up{\Upsilon}
\def\vp{\varphi}
\def\ve{\varepsilon}
\def\wh{\widehat}
\def\br{\breve}
\def\wt{\widetilde}
\def\ov{\overline}

\def\BC{{\mathbb C}}
\def\BR{{\mathbb R}}
\def\BN{{\mathbb N}}

\def\clp{{\mathcal P}}

\def\cln{{\mathcal N}}

\def\E{\mathrm{e}}
\newtheorem{Pa}{Paper}[section]
\newtheorem{Tm}[Pa]{{\bf Theorem}}

\newtheorem{Cy}[Pa]{{\bf Corollary}}
\newtheorem{Rk}[Pa]{{\bf Remark}}

\newtheorem{Dn}[Pa]{{\bf Definition}}
\newtheorem{Pn}[Pa]{{\bf Proposition}}

\title{Weyl functions of Dirac systems and of their generalizations: integral representation,
inverse problem, and discrete interpolation}

\author{B. Fritzsche, B. Kirstein, A.L. Sakhnovich}

\date{}
\begin{document}
\maketitle

\begin{abstract} 
Self-adjoint Dirac systems and  subclasses of canonical systems,
which generalize Dirac systems are studied.
Explicit and global solutions of direct and inverse problems
are obtained. A local Borg-Marchenko-type theorem,
integral representation of the Weyl function, and results
on interpolation of  Weyl functions are also derived.
\end{abstract}



\section{Introduction} \label{intro}
\setcounter{equation}{0}
Weyl functions (also called Weyl-Titchmarsh or $M$-functions) and their generalizations are
an important and much used  tool in the spectral theory of differential equations
(see a necessarily small part of recent papers and books on this topic and various references therein 
 \cite{AP, AD, AMR, CG, CGZ, EK, FKS, Ges, GeSi, GZ, MST, SaA1, SaAn, SaL4, T}).    
Following the seminal work \cite{Kr} and  more general constructions  in \cite{SaL2,
SaL4} (see also some references therein), one can use structured operators to solve inverse
problems for Krein, Dirac, canonical, non-classical and non-self-adjoint systems and their
discrete analogues (see, for instance, \cite{AGKLS, FKRS, FKS, SaA1, SaA2, SaAn, SaAerd}).
It was proved in \cite{FKS, SaA0, SaA1, SaA2, SaAn, SaAerd} that the kernels of these
structured operators are connected with the Weyl functions via some kinds of Fourier
transformations and can be recovered directly from the Weyl functions.
It is also true that the Weyl functions can be represented as Fourier transformations
of the corresponding kernels of the structured operators.
Important related works  \cite{GeSi, Si} on the $A$-amplitudes for  Schr\"odinger 
operators gave rise to a whole series of interesting papers and results
on the high energy asymptotics of Weyl functions and  local Borg-Marchenko-type
uniqueness theorems (see, for instance, \cite{CG0, CGZ, CGZ2, LaWo} and references therein). See also
further  discussion of $A$-amplitudes in \cite{AMR}.
Finally,  the interpolation 
of Weyl functions for scalar Schr\"odinger operators using the values of Weyl functions on a discrete
countable set was dealt with in two interesting papers \cite{BT, RT}.
The second and more general paper \cite{RT} is based on the $A$-amplitude representation
of Weyl functions.

 In this paper we consider self-adjoint Dirac systems and  subclasses of canonical systems,
which generalize Dirac systems. We consider direct and  inverse problems, integral representations
of Weyl functions, and interpolation of Weyl functions.
Connections with Schur coefficients and $A$-amplitude are discussed.

In Section \ref{saDir} we formulate some previous results for Dirac systems, 
prove formula \eqref{d3}, which solves the inverse problem in a more
general setting than previously considered, and derive for Dirac systems analogs of some 
$A$-amplitude-type results for Schr\"odinger equations.
In Section \ref{Can} subclasses of canonical systems are studied and
the case of rational Weyl matrix functions is discussed.
Explicit solutions of direct and inverse problems are given for that case.

A general case, where Weyl functions are not necessarily rational,
is considered in Section \ref{CanG}.
An integral representation
of the Weyl function, a  solution of the inverse problem for generalizations
of Dirac systems directly via Weyl functions,
and a local Borg-Marchenko-type theorem are obtained 
there.
Section \ref{Intrp} is dedicated to interpolation of Weyl functions.

As usual, by $\BC$ we denote the complex plane,  by $\BC_+$ we denote the open upper
semi-plane, $\BR$ is the real axis, $\BR_+$ is the positive real axis,
$\BN$ is the set of positive integers, 
and $\BN_0$ is the set of non-negative integers.  We have $\b^{\prime}=\frac{d}{dx}\b$.
By $L^2_{p\times p}(-\infty,\, \infty)$ we mean the space of $p \times p$
matrix functions with entries which belong to $L^2(-\infty,\, \infty)$. The $p\times p$ identity
matrix is denoted by $I_p$,  the identity operator is denoted by $I$, and diag
denotes a diagonal or block diagonal matrix.
The set of bounded linear operators acting
from $H_1$ to $H_2$ is denoted by $\{H_1, H_2\}$  and $\s$ is used to denote the spectrum.

\section{Self-adjoint Dirac system: \\ $A$-amplitude and inverse problem} \label{saDir}
\setcounter{equation}{0}
Consider the matrix self-adjoint   Dirac system
\begin{equation}       \label{1.1}
\frac{d}{dx}u(x, z)=i\big(z j+jV(x)\big)u(x, z), \quad x\geq 0,
\end{equation}
where
\begin{equation}   \label{1.2}
j =
\left[
\begin{array}{cc}
I_{p} & 0 \\ 0 & -I_{p}
\end{array}
\right], \hspace{1em}
V= \left[\begin{array}{cc}
0&v\\v^{*}&0\end{array}\right],
 \end{equation}
$I_{p}$
is the $p \times p$ identity matrix, $V$ is an
 $m \times m $ ($m=2p$) matrix function with the
 $p \times p$ block entry $v$ (often referred to as the potential)
 being locally summable.  If $v$ is summable on each finite interval
 $[0, \, l]$, that is, if $v$ is locally summable, there is (see Proposition 5.2 and its proof in \cite{SaA2}) a unique $p \times p$ matrix function 
 $\vp(z)$ such that
 \begin{equation} \label{1.3}
\int_0^\infty
\left[ \begin{array}{lr}
I_p & i \varphi (z)^* \end{array} \right]
  \widehat u(x, z)^*
\widehat u(x, z)
 \left[ \begin{array}{c}
I_p \\ - i \varphi (z) \end{array} \right]
dx            < \infty
 \quad (z \in {\BC}_+),
\end{equation}
where $\widehat u$ is the fundamental solution of (\ref{1.1})
normalized by the  condition 
\begin{align}\label{1.4}
\widehat u(0, z)=K^*, \quad K:=   \frac{1}{\sqrt{2}}       \left[
\begin{array}{cc} I_p &
-I_{p} \\ I_{p} & I_p
\end{array}
\right].
\end{align}
\begin{Dn} The $p \times p$ matrix function $\vp(z)$ $\,(z\in \BC_+)$, which satisfies \eqref{1.3},  is called the Weyl function
of system \eqref{1.1}. 
\end{Dn}
The Weyl function  of  \eqref{1.1} is holomorphic in $\BC_+$ and $\Im \vp \geq 0$, which we discuss in detail in Section \ref{Intrp}.

We always assume that $v$ is measurable.
When $v$ is also locally bounded, one can recover it from the spectral function or directly
from the Weyl function using Theorems 4.2 or 5.4 from \cite{SaA2}, respectively.  
To recover $v$ from the Weyl function one should notice that  the Weyl function
of system \eqref{1.1} on the semi-axis is a so called Weyl point
and, according
to Proposition 5.2 and formula (5.3) from \cite{SaA2},
Theorem 5.4 from \cite{SaA2}
is applicable. 

Let us describe the procedure
suggested in Theorem 5.4 \cite{SaA2} to recover system on  any intervals.
For this purpose we need a family of operators 
which are associated with system \eqref{1.1} and have difference kernels:
\[
S_l=\frac{d}{dx}\int_0^l s(x-t) \, \cdot \, dt, \quad s(x)=-s(-x)^*,
\]
where
\begin{equation} \label{1.5}
s(x)= \left( \frac{d}{dx} \frac{i}{4 \pi}
e^{ \eta x}
{\mathrm{l.i.m.}}_{a\to \infty} \int_{- a}^{a}e^{-i \zeta x}
(\zeta +i \eta)^{-2} \varphi (\zeta +i \eta)d \zeta \right)^*, \quad \eta>0.
\end{equation}
Here $(\zeta +i \eta)^{-2} \varphi (\zeta +i \eta)\in L^2_{p\times p}(-\infty,\, \infty)$
for every fixed $\eta>0$, l.i.m. denotes the  entrywise limit in the norm of  $L^2(0,\, l)$
$(l>0)$, and  l.i.m. in the right-hand side of \eqref{1.5} is differentiable.
Moreover,  $s(x)$ is boundedly differentiable, $s(0)=\frac{1}{2}I_p$, and $S_l>0$
for all $l>0$. Thus, we have
\begin{align}\label{1.6}&
S_l=I+\int_0^lk(x-t) \, \cdot \, dt>0, \quad
k(x):=s^{\prime}(x) \,\, \big(s^{\prime}=\frac{d}{dx}s\big), \quad k(x)=k(-x)^*, 
\end{align}
and $S_l, \, S_l^{-1}\in \{L^2_p(0,l), L^2_p(0,l)\}$. Notice, that for convenience purposes
$\vp$, $s$, and  $S$  in the present
paper differ by scalar constant factors from the objects denoted by the same 
letters in  \cite{SaA2}.

We now introduce the matrix functions 
 \begin{align}\label{1.6'}&
\t_1(x)= [I_p \quad 0]\wh u(x,0), \quad  \t_2(x)= [0 \quad I_p]\wh u(x,0).
\end{align}
According to formula (4.16) from \cite{SaA2} we have
 \begin{align} \label{1.7}&
\t_2(x)= \frac{1}{\sqrt{2}} \left( [-I_p \quad I_p]
- \int_0^{2x} k(t)^* S_{2x}^{-1}[  2s(t) \quad I_p]dt \right),
 \end{align}
where $S_{r}^{-1}$ $( r=2x)$ is applied to $[  2s(t) \quad I_p]$ columnwise.
It easily follows from \eqref{1.1} and \eqref{1.4} that
\begin{equation} \label{1.8}
\wh u(x,0)^*j\wh u(x,0)\equiv J, \quad \wh u(x,0)J\wh u(x,0)^*\equiv j, \quad
J:=\left[
\begin{array}{cc}
0& I_{p}  \\ I_{p} & 0
\end{array}
\right].
\end{equation}
Using \eqref{1.1}, \eqref{1.6'}, and the second relation in \eqref{1.8} we obtain
 \begin{align} \label{1.9}&
 v(x)=i\t_1^{\prime}(x)J\t_2(x)^*.
 \end{align}
 Finally, in view of  \eqref{1.1}, \eqref{1.4}, \eqref{1.6'}, and the second relation in \eqref{1.8}
we get equalities
 \begin{align} \label{1.10}&
 \t_1(0)= \frac{1}{\sqrt{2}}  [I_p \quad I_p], \quad \t_1(x)J\t_2(x)^*\equiv 0, \quad \t_1^{\prime}(x)J\t_1(x)^*\equiv 0,
 \end{align} 
 which uniquely determine $\t_1$ assuming that $\t_2$ is already given.
Rewrite Theorem 5.4 from \cite{SaA2} in a slightly modified form. 
\begin{Tm} \label{saDip}\cite{SaA2}  Let $\vp$ be the Weyl function of the matrix Dirac system \eqref{1.1}, where  the 
potential $v$ is locally bounded, that is,
 \begin{align} \label{1.10'}&
\sup_{0< x<l}\|v(x)\|<\infty \quad {\mathrm{for \,\, any}} \quad l>0.
 \end{align} 
Then $v$ is uniquely recovered from $\vp$
by the formulas \eqref{1.9}, \eqref{1.7}, and \eqref{1.10}, whereas $s$ and $S_{2x}$ in \eqref{1.7}  are given by
\eqref{1.5} and \eqref{1.6}, respectively.
 \end{Tm}

\begin{Rk}\label{ampl}
Notice that according to formula (3.26) from \cite{SaA2} we have
\begin{align} \label{1.11}&
\vp(z)=2z\int_0^{\infty}e^{izx}s(x)^*dx,
 \end{align} 
 that is, the kernel $s$, which is used in Theorem \ref{saDip}  to solve the inverse
problem, is a Dirac system's analog of the $A$-amplitude. 
(To be more precise, $k=s^{\prime}$ is the analog.)
Following M.G. Krein, the matrix functon $k$  is called an accelerant in \cite{AGKLS}. 
\end{Rk}
Clearly,
formulas \eqref{1.5} and \eqref{1.11} are closely related (and \eqref{1.11} was used
in \cite{SaA2} to derive \eqref{1.5}.)  Another closely related formula is representation
(5.5) from \cite{SaA2}:
\begin{equation} \label{1.12}
\vp(z)=z^2 \int_0^{\infty}e^{izx}\chi(x)dx, \quad e^{-\eta x}\chi(x)\in L^2_{p \times p}(0, \, \infty) \,\, {\mathrm{for \,\, all}}\,\, \eta>0,
\end{equation}
 where
\begin{align} \label{1.13}&
\chi(x):=-2i\int_0^xs(t)^*dt.
 \end{align} 

According to p. 341 in \cite{SaA2} the operators $S_l$ admit a triangular factorisation
\begin{align} \label{d1}&
S_l=E_l^{-1}\big(E_l^{-1}\big)^*, \quad E_l=I+\int_0^x E(x,t) \, \cdot \, dt, 
\\ \label{d2} &
  E_l^{-1}=I+\int_0^x \Ga (x,t) \, \cdot \, dt, \quad \sup_{0<t<x<l}(\|E(x,t)\|+\|\Ga(x,t)\|)<\infty, 
 \\ \label{d2'} & 
 E_l^{-1}\t_{1}(x/2)=\frac{1}{\sqrt{2}}[2s(x) \quad I_p]  \quad (0<x<l) ,
 \end{align} 
where $E_l,\, E_l^{-1} \in \{L^2_p(0,l), L^2_p(0,l)\}$
and  $E_l^{-1}$ is applied to $\t_{1}(x/2)$ columnwise.
 (Note that the expression for the triangular factor $V$ of  $S_l^{-1}$ on p. 341 in \cite{SaA2}
is slightly more complicated but according to \cite{SaL0}
we can put $\kappa(x)\equiv 1/\sqrt{2}$ in that expression.
Taking into account also that  $S_l$ is here two times less
than $S_l$  in \cite{SaA2} and so $E_l$ is here $\sqrt{2}$ times larger than its analog
$V$ in \cite{SaA2}, we get precisely  \eqref{d1}.)   
By \eqref{d1} and \eqref{d2}  we have
\begin{align} \label{d2!}&
S_l^{-1}=E_l^{*}E_l, \quad l>0,
 \end{align} 
and the matrix function $E$ is locally bounded, that is, we can use equality \eqref{d2!} to determine
pointwise the kernel of the integral operator in the left hand-side.  Specifically, we get
\begin{align} \nonumber
\big(S_l^{-1}f\big)(x)=&f(x)+\int_0^x\Big(E(x,t)+\int_x^l E(r,x)^*E(r,t)dr\Big)f(t)dt \\
 \label{d2!!}&
+
\int_x^l\Big(E(t,x)^*+\int_t^l E(r,x)^*E(r,t)dr\Big)f(t)dt.
 \end{align}

Now,
the next result  follows from Theorem \ref{saDip} and formula \eqref{d2!!}.
\begin{Tm} \label{saDnf}  Let $\vp$ be the Weyl function of the matrix Dirac system \eqref{1.1}, where  the 
potential $v$ satisfies \eqref{1.10'}.
Then $v$ is uniquely recovered from $\vp$
by the formula
\begin{align} \label{d3}&
v\big(\frac{l}{2}\big)=2i\big(S_l^{-1}k\big)(l), 
 \end{align} 
where $k$ and $S_l$ are given by
\eqref{1.5} and \eqref{1.6}.
 \end{Tm}
\begin{proof}.  By \eqref{d2!} we get  $S_r^{-1}=E_r^*E_r$. 
Note also that $(E_rf)(t)=(E_{l}f)(t)$ for $0\leq t\leq r\leq l$.
Hence, formulas
\eqref{1.7}  and \eqref{d2'}  imply
 \begin{align} \label{d4}&
\t_2^{\prime}(x)=-\sqrt{2}\Big(E_l k\Big)(2x)^*
\Big(E_l[2s \quad I_p]\Big)(2x)=-2\Big(E_l k\Big)(2x)^*\t_1(x), \quad 2x<l.
 \end{align} 
In view of \eqref{1.9} and the second equality in \eqref{1.10} one can see that
\begin{align} \label{d5}&
\t_2^{\prime}J\t_1^*=-\t_2J(\t_1^{\prime})^*=-iv^*.
 \end{align} 
By the first equality in  \eqref{1.6'} and the second equality in \eqref{1.8} we have 
\begin{align} \label{d6}&
\t_1J\t_1^*\equiv I_p.
 \end{align} 
From \eqref{d4}-\eqref{d6} we obtain $v\big({x}\big)=2i\Big(E_l k\Big)(2x)$ or, equivalently,
\begin{align} \label{d7}&
v\big(\frac{x}{2}\big)=2i\Big(E_l k\Big)(x), \quad 0<x<l.
 \end{align} 
 In other words, for almost all values of $x$ we have
 \begin{align} \label{d7'}&
v\big(\frac{x}{2}\big)=2i\Big(k(x)+\int_0^xE(x,t)k(t)dt\Big).
 \end{align} 
Putting $x=l$ in \eqref{d2!!} and  \eqref{d7'}  we see that \eqref{d3} is true.
\end{proof}
Formulas related to \eqref{d3} can be found in \cite{AGKLS, FKS}
for the cases of continuous $v$ and of a skew-self-adjoint Dirac system,
respectively.

The case of a continuous $v$ was treated in \cite{AGKLS}.  If $v$ is continuous, then
by Theorem 1.2 in \cite{AGKLS}, the matrix function $k$ is also continuous with
a possible jump discontinuity at the origin. Moreover, it follows from the proof 
of Theorem 1.2 in \cite{AGKLS} that the kernels  $E(x,t)$ and $\Ga(x,t)$ of the
operators $E_l$ and $E_l^{-1}$, respectively, are  continuous in the triangles 
$0 \leq t \leq x\leq l$.  In particular, $E(x,0)$ is well-defined. 
According to formulas (1.13) and (2.12) from \cite{AGKLS} we have
\begin{align} \label{d9}&
v(x)=-2iE(2x,0).
 \end{align} 
Moreover, since $E(x,t)$ is continuous, formula \eqref{d7}
 holds pointwise. By applying $E_l^{-1}$ to both sides of \eqref{d7}
 we get the following analog of  the representation 
  of $A$-amplitude of Schr\"odinger operator  in \cite{AMR}.  (See formula (3.11) in \cite{AMR},
which is essential in the  proof of Theorem 1 in the same paper.)
  \begin{Cy} Let  the potential $v$ of system \eqref{1.1} be continuous.
  Then the matrix function $k$ given by  \eqref{1.5} and \eqref{1.6}
  admits the representation
  \begin{align}\nn
  k(2x)&=-\frac{i}{2}\Big(v(x)+2\int_0^x\Ga\big(2x,2t\big)v(t)dt\Big)
 \\  \label{d10}&=-\frac{i}{2}\Big(v(x)+2\int_0^x\Ga\big(2x,2(x-t)\big)v(x-t)dt\Big).
  \end{align}
  The matrix functions $k$ and $\Ga$ are continuous and \eqref{d10} holds for all
  $x\in [0,\, \infty)$.
  \end{Cy}



\section{Canonical system. Explicit formulas} \label{Can}
\setcounter{equation}{0}
\subsection{Canonical system and Schur coefficients}
In this section we consider a canonical system
\begin{equation}\label{v1.3}
\frac{d}{d x}w(x, z   )=i z   JH(x)w(x, z   ),
\quad
H(x)=\b(x)^*\b(x), \quad   x \geq 0,
\end{equation}
where $J$ is given by \eqref{1.8}, $\b$ is  an $p \times m$ ($m=2p$) matrix function
such that
\begin{align}\label{c1}&
\b(x)J\b(x)^*=D=\diag\{d_1, \, d_2, \ldots , d_p\},
\end{align}
and diag denotes a diagonal matrix.
Weyl  function of the canonical system on the semi-axis $x \geq 0$ is a $p \times p$ 
holomorphic  matrix function $\vp( z)  $,
which satisfies the condition \cite{SaL4}
\begin{equation} \label{v1.6}
\int_0^\infty
\left[ \begin{array}{lr}
I_p &i \varphi ( z)  ^* \end{array} \right]
  w(x,  z)  ^*H(x)w(x,  z)  
\left[ \begin{array}{c}
I_p \\ -i \varphi ( z)   \end{array} \right]
dx            < \infty, \quad z\in \BC_+,
\end{equation}
where $w$ is the fundamental solution of  \eqref{v1.3} normalized by $w(0,z)=I_{2p}$.
The Dirac system of the previous section, is equivalent  \cite{SaA2}
to  canonical system
\begin{equation}\label{c2}
\frac{d}{d x}\wh w(x, z   )=i z   J \wh H(x) \wh w(x, z   ),
\quad
\wh H(x)=\wh u(x,0)^*\wh u(x,0), \quad   x \geq 0.
\end{equation}
Taking into account \eqref{1.8} one can see that
a simple multiplication $w(x,z)=e^{izx}\wh w(x,z)$ of solution $\wh w$ by scalar factor
transforms system \eqref{c2} into system \eqref{v1.3}, where
\begin{align}\label{c3}&
D=2I_p, \quad \b(x)=\sqrt{2}\t_1(x),
\end{align}
and $\t_1$ is given by \eqref{1.6'}. Nevertheless, additional Weyl functions appear
as a result of this transformation and the Weyl  theory of system \eqref{v1.3} requires
some separate study even in this particular case.

Before we turn to the general case \eqref{c1}, we consider the canonical system \eqref{v1.3}, \eqref{c3} generated by Dirac system in greater detail.
Separating  $\b$ into two $p \times p$ blocks, namely $\b_1$ and $\b_2$,
and using \eqref{c1} and  \eqref{c3}  we obtain
\begin{align}\label{r17}&
\det \, \b_2(x) \not =0, \quad \Re\big(\b_2(x)^{-1}\b_1(x)\big)>0.
\end{align}
By \eqref{1.1}, \eqref{1.6'}, and  \eqref{1.8}  we get
\begin{align}\label{r16}
&\b^{\prime}(x)J\b(x)^*=2\t_1^{\prime}(x)J\t_1(x)^*=2iv(x)[0 \quad I_p]j\left[
\begin{array}{c}
I_{p}  \\ 0 
\end{array}
\right]=0.
\end{align}
According to  \eqref{r16} we have $\b_1^{\prime}\b_2^*+\b_2^{\prime}\b_1^*=0$, that is,
\begin{align}
\label{r18} &
\b_1^{\prime}(x)=-\b_2^{\prime}(x)\b_1(x)^*\big(\b_2(x)^*\big)^{-1}.
\end{align}
It follows from \eqref{r18}    that 
$\big(\b_2^{-1}\b_1\big)^{\prime}=-\b_2^{-1}\b_2^{\prime}\big(\b_2^{-1}\b_1+\b_1^*(\b_2^*)^{-1}\big)$.
Moreover, formulas \eqref{1.4}, \eqref{1.6'}, and \eqref{c3} imply $\b_2(0)=I_p$.
Hence, we can recover $\b_2$ from $\b_2^{-1}\b_1$ using the following differential equation
and initial condition
\begin{align}\label{r19}
&\b_2^{\prime}=-\b_2\big(\b_2^{-1}\b_1\big)^{\prime}\big(\b_2^{-1}\b_1+\b_1^*
(\b_2^*)^{-1}\big)^{-1}, \quad \b_2(0)=I_p.
\end{align}
\begin{Rk}\label{Schur}  Note that  $\b_2^{-1}\b_1$ is the continuous analog
of the (well-known in the discrete case) Schur coefficients. Similar to the case of
Toeplitz matrices, the function $\b_2^{-1}\b_1$  has the property $\Re\big(\b_2^{-1}\b_1\big)>0$.
Moreover, system  \eqref{v1.3} corresponding to operators with difference kernels
is uniquely recovered from $\b_2^{-1}\b_1$ via formulas \eqref{r19} and
$\b_1=\b_2\big(\b_2^{-1}\b_1\big)$.
\end{Rk}

\subsection{Explicit formulas}

In the case of the rational Weyl functions and corresponding Hamiltonians $H$ one can obtain explicit 
formulas for solutions of the direct and inverse problems
using a GBDT version of the B\"acklund-Darboux transformation \cite{SaA2c, SaA3}.
See \cite{Ci,D, GKS6, GeT, Gu, MS, SaA2c, ZM} and references therein for various versions
of the  B\"acklund-Darboux transformation and commutation methods. For explicit formulas
via inversion of semiseparable operators see, for instance, \cite{AGKLS, FKS, ASAK}. 
Here we shall apply the procedure from \cite{SaA2c} to the
initial system
\begin{equation} \label{v1.8}
\frac{d}{d x}w_0(x, z   )=i z   JH_0w_0(x, z   ),
\quad
H_0 \equiv \left[ \begin{array}{c}  D/2 \\
I_p 
\end{array} \right]\left[ \begin{array}{lr}  D/2 &
I_p 
\end{array} \right] \quad (\det D \not=0),
\end{equation}
where $D=D^*=\diag\{d_1, \, d_2, \ldots , d_p\}$ is a fixed diagonal matrix.
(Notice that $J$ in (\ref{v1.3})  and  in (\ref{v1.8})  is slightly different from $J$ in  \cite{SaA2c}.)

We consider a fixed integer $n>0$, an $n \times n $ matrix $\a$ and two $n \times p$ matrices
$\Lam_1(0)$ and  $\Lam_2(0)$, such that
\begin{equation} \label{v1.9}
\a - \a^*=i\Lam(0)J\Lam(0)^*, \quad \Lam(0):=[\Lam_1(0) \quad \Lam_2(0)].
\end{equation}
We now introduce matrix functions $\Lam(x)$ and $\S(x)$ with the equalities
\begin{equation} \label{v1.10}
\frac{d}{d x}\Lam=-i \alpha \Lam J H_0, \quad \Lam(0)=[\Lam_1(0) \quad \Lam_2(0)];  
\end{equation}
\begin{equation} \label{v1.11}
  \frac{d}{d x}\S=\Psi_1\Psi_1^*, \quad \S(0)=I_n, 
\quad \Psi_k(x):=\Lam(x)\left[ \begin{array}{c}I_p\\ (- 1)^{k+1} D/2 
\end{array} \right], \quad k=1,2.
\end{equation}
The first two equalities in (\ref{v1.11})  imply $\S(x)\geq I_n$. Thus, $\S$ is invertible and
we can put
\begin{equation} \label{v1.11.1}
H(x)=v_0(x)^*H_0v_0(x), 
\end{equation}
where
\begin{eqnarray} \label{v1.11.2}&&
  \frac{d}{d x}v_0=- q_0v_0, \quad  v_0(0)Jv_0(0)^*=J,
\\  \label{v1.11.3} &&
 q_0(x):=J\Lam(x)^*\S(x)^{-1}\Lam(x)JH_0-
  JH_0J\Lam(x)^*\S(x)^{-1}\Lam(x).
\end{eqnarray}
The transfer matrix function, which we use in GBDT, has the form
\begin{equation} \label{v1.11.4}
w_{\a}(x, z)  =I_m-iJ\Lam(x)^*\S(x)^{-1}(\a - z I_n)^{-1}\Lam(x).
\end{equation}
In view of  the second equality in (\ref{v1.8}), the first equality in (\ref{v1.10})
and the third equality in (\ref{v1.11}), we see that  
$\Big(\frac{d}{dx} \Lam \Big)J \Lam^*=-i \a \Psi_1 \Psi_1^*$.  Therefore,
taking into account  (\ref{v1.9}) and (\ref{v1.11}) we get
\[
\frac{d}{d x}\Big(\a \S - \S \a^* \Big)=i\frac{d}{d x}\Big(\Lam J \Lam^* \Big), \quad
\a \S(0) - \S(0) \a^*=i\Lam(0) J \Lam(0)^*.
\]
The two equalities above
imply 
\begin{equation} \label{vd2}
\a \S(x) - \S(x) \a^*=i\Lam(x) J \Lam(x)^*.
\end{equation}
By direct calculation\cite{SaL2, SaL4}, it follows from (\ref{v1.11.4}) and (\ref{vd2})
that
\begin{equation} \label{vd3}
w_{\a}(x,\ov  z)  ^*Jw_{\a}(x, z)  =J.
\end{equation}
It follows from (\ref{v1.3}) and (\ref{v1.11.2}), respectively, that
\begin{equation} \label{vd4}
w(x,\ov  z)  ^*Jw(x, z)  =J, \quad  v(x)^*Jv(x)=J.
\end{equation}
According to Theorem 1  in \cite{SaA2c} the fundamental solution $w$
of  the canonical system given by (\ref{v1.3}) and (\ref{v1.11.1})-(\ref{v1.11.3})   is expressed
via $w_{\a}$:
\begin{equation} \label{v1.11.5}
w(x, z)  =v_0(x)^{-1}w_{\a}(x, z)  w_0(x, z)  w_{\a}(0, z)  ^{-1}v_0(0),
\end{equation}
where $w_0$ is the fundamental solution of the initial system (\ref{v1.8}).
Here the fundamental solutions are normalized by the initial condition
$w(0, z)  =w_0(0, z)  =I_m$ $\, (m=2p)$.

\begin{Rk} \label{ReExpl} If $\det \a \not= 0$, by formula (8)  in \cite{SaA2c} and by formula (\ref{v1.11.2}) we have
\begin{equation} \label{v1.11.6}
v_0(x)=w_{\a}(x,0)M, \quad M:=v_0(0), \quad MJM^*=J.
\end{equation}
Moreover, for the case of the initial system (\ref{v1.8}) the matrix functions $w_0$, $\Lam$ and $\S$
can be constructed explicitly. Hence, assuming $\det \a \not= 0$, one gets explicit expressions
for $H(x)$ and $w(x, z)  $ via formulas (\ref{v1.11.1}), (\ref{v1.11.4}), (\ref{v1.11.5}),  and (\ref{v1.11.6}).

Indeed, put
\begin{equation} \label{v1.11.7}
Z= \left[ \begin{array}{lr} I_p & I_p \\
D/2 & -D/2
\end{array} \right], \quad \br D=\left[ \begin{array}{lr}  D & 0
\\0 & 0
\end{array} \right].
\end{equation}
It follows that
\begin{equation} \label{v1.11.8}
Z^{-1}= \diag\{D^{-1},\, D^{-1}\} \left[ \begin{array}{lr} D/2 &I_p \\
D/2 & -I_p
\end{array} \right].
\end{equation}
In view of  (\ref{v1.8}), (\ref{v1.11.7}),  (\ref{v1.11.8}) and normalization condition $w_0(0, z)  =I_m$,
it is immediately apparent that
\begin{equation} \label{v1.11.9}
w_0(x, z)  =Ze^{i z x \br D}Z^{-1}.
\end{equation}
According to  (\ref{v1.10}), to the second equality in (\ref{v1.8}) and to the third equality in (\ref{v1.11}), we have
\begin{equation} \label{v1.12}
\frac{d}{d x}\Psi_1=-i \alpha \Psi_1 D, \quad  \Psi_1(0) =  \Lam_1(0) +\frac{1}{2} \Lam_2(0) D,
\quad \Psi_2(x) \equiv  \Lam_1(0) - \frac{1}{2} \Lam_2(0) D.
\end{equation}
The  expression for $\Psi_1$ follows directly from  (\ref{v1.12}):
\begin{equation} \label{v1.13}
\Psi_1(x)=[ \exp(-i d_1 x \a) f_1\quad  \exp(-i d_2 x \a)f_2  \quad \ldots \quad  \exp(-i d_p x \a)f_p],
\end{equation}
\begin{equation}
 \label{v1.14}
[  f_1\quad f_2  \quad \ldots \quad  f_p]:= \Psi_1(0) .
\end{equation}
Formulas  (\ref{v1.12}) and (\ref{v1.13}) provide explicit expressions for $\Psi_1$ and $\Psi_2$.
Finally, from  (\ref{v1.11}) and (\ref{v1.11.7}) we obtain
\begin{equation} \label{v1.15}
\Lam(x)=[\Psi_1(x) \quad \Psi_2(x)]Z^{-1}, \quad \S(x)=I_n+\int_0^x\Psi_1(t)\Psi_1(t)^*dt.
\end{equation}
\end{Rk}
\begin{Dn}\label{Dn1}
A canonical system given by (\ref{v1.3}) and (\ref{v1.11.1}), where $H_0$
is defined in (\ref{v1.8}), 
$v_0$ is defined by (\ref{v1.11.2}) and (\ref{v1.11.3}), and
$\Lam (x)$ and $\S(x)$ in (\ref{v1.11.3}) are defined by  (\ref{v1.10}) and (\ref{v1.11}) or, equivalently,
by \eqref{v1.12}-\eqref{v1.15},
is said to be determined by the parameter matrices $\a$ and $\Lam (0)$ such that
 (\ref{v1.9}) holds.
\end{Dn}
Further we always assume that  (\ref{v1.9}) is valid.

\subsection{Direct problem: explicit solutions}
Remark \ref{ReExpl} leads us to our next proposition.
\begin{Pn}  \label{PnWe} Let $v_0(0)=I_m$ $(m=2p)$. Then the matrix functions
\begin{align} \label{v1.16}&
\vp( z)  =-\frac{i}{2}D+\Psi_1(0)^*(\g - z I_n)^{-1}\Psi_2, \\
\label{c1.16}&
 \wh \vp( z)  =\frac{i}{2}|D|+\wh \Psi_1(0)^*(\wh\g - z I_n)^{-1}\wh \Psi_2,
\end{align}
where 
\begin{align} & \label{v1.17}
\g=:\a-i\Psi_2\Lam_2(0)^*, \quad \Psi_1(0)= \Lam_1(0) + \frac{1}{2} \Lam_2(0) D,\\
& \label{v1.17'}
 \Psi_2= \Lam_1(0) - \frac{1}{2} \Lam_2(0) D \equiv  \Lam_1(x) - \frac{1}{2} \Lam_2(x) D,\\
\label{c1.17}&
\wh \g=:\a-i\wh\Psi_2\Lam_2(0)^*, \quad \wh\Psi_1(0)= \Lam_1(0) - \frac{1}{2} \Lam_2(0) |D|,\\
& \label{c1.17'}
\wh \Psi_2= \Lam_1(0) + \frac{1}{2} \Lam_2(0) |D|, \quad |D|=\diag\{|d_1|, |d_2|, \ldots, |d_p| \}
\end{align}
are Weyl functions of  the canonical system determined by $\a$ and $\Lam (0)$.

If $D<0$, then the Weyl functions $\vp( z)  $ and $\wh \vp( z)  $
coincide and, moreover,  the  Weyl function is unique. 

The matrices $\g$ and $\wh \g$ given by   (\ref{v1.17})  and  (\ref{c1.17}), respectively, satisfy matrix identities
\begin{equation} \label{v1.17''}
\g - \g^*=i\Lam_2(0)D\Lam_2(0)^*, \quad \wh \g^* -\wh  \g=i\big(\wh\Psi_2- \wh\Psi_1(0)\big)|D|^{-1}\big(\wh\Psi_2- \wh\Psi_1(0)\big)^*,
\end{equation}
that is, for $D<0$ we obtain $\s(\g) \cap \BC_+=\emptyset$ and $\vp$ 
given by \eqref{v1.16} does not have
singularities in $\BC_+$. 

\end{Pn}
\begin{proof}. The proof that $\vp$ 
is a Weyl function is somewhat similar to the proof of the corresponding fact for the
self-adjoint Dirac  system \cite{GKS1}. Put
\begin{equation} \label{v1.18}
 \Om( z)  :=w_{\a}(0,  z)  Z, \quad
 \Om_2( z)  =\left[ \begin{array}{c}   \om_1( z)   \\   \om_2( z)  
\end{array} \right]:=
\Om( z)  
\left[ \begin{array}{c}  0 \\  I_p
\end{array} \right],
\end{equation}
where $  \om_k$ are $p \times p$ blocks of $   \Om_2$. In view of 
(\ref{v1.11.5}), (\ref{v1.11.9})  and  (\ref{v1.18}) we have
\begin{equation} \label{v1.19}
w(x, z)   \Om_2( z)  =v_0(x)^{-1}w_{\a}(x, z)  Z\left[ \begin{array}{c}  0 \\  I_p
\end{array} \right]=v_0(x)^{-1}w_{\a}(x, z)  \left[ \begin{array}{c}  I_p \\
-D/2
\end{array} \right].
\end{equation}
Notice also that by  (\ref{v1.11}) and  (\ref{v1.11.4}) the equality
\begin{equation} \label{v1.20}Q(x):=\left[ \begin{array}{lr}  D/2 &
I_p 
\end{array} \right]w_{\a}(x, z)  \left[ \begin{array}{c}  I_p \\
-D/2
\end{array} \right]=-i\Psi_1(x)^*\S(x)^{-1}(\a - z I_n)^{-1}\Psi_2
\end{equation}
holds.
From  (\ref{v1.8}), (\ref{v1.11.1}),  (\ref{v1.19}), and  (\ref{v1.20})  it follows that
\begin{eqnarray} \label{v1.21}
&&\Om_2( z)  ^*w(x, z)  ^*H(x)w(x, z)   \Om_2( z)  = Q(x)^* Q(x)
\\ \nonumber &&
=\Psi_2^*(\a^* - 
\ov z I_n)^{-1}\S(x)^{-1}\Psi_1(x)\Psi_1(x)^*\S(x)^{-1}(\a - z I_n)^{-1}\Psi_2
.
\end{eqnarray}
According to (\ref{v1.11}) we have 
\begin{align}\label{c4}
\frac{d}{d x}\S^{-1}=-\S^{-1}\Psi_1\Psi_1^*\S^{-1}, \quad \S(0)=I_n, \quad \S(x)\geq I_n.
\end{align}
Hence, from  (\ref{v1.21}) it follows that
\begin{eqnarray} \label{v1.22}&&
\int_0^r \Om_2( z)  ^*w(x, z)  ^*H(x)w(x, z)   \Om_2( z)  dx =\Psi_2^*(\a^* - 
\ov z I_n)^{-1}
\\ \nonumber &&
\times \big(\S(0)^{-1}-\S(r)^{-1}\big)(\a - z I_n)^{-1}\Psi_2
\leq \Psi_2^*(\a^* - 
\ov z I_n)^{-1}(\a - z I_n)^{-1}\Psi_2.
\end{eqnarray}
Compare (\ref{v1.6}) and  (\ref{v1.22}) to see that the function
\begin{equation} \label{v1.23}
\vp( z)  =i \om_2( z)   \om_1( z)  ^{-1},
\end{equation}
where $ \om_k$ are the blocks of $ \Om_2$, 
satisfies (\ref{v1.6}) (excluding, possibly, a finite number of points), i.e., $\vp$ of the form \eqref{v1.23}
is a Weyl function of system (\ref{v1.3}),  (\ref{v1.11.1}).

Next, let us show that the right-hand sides of (\ref{v1.16}) and  (\ref{v1.23})
coincide. Since $v_0(0)=I_m$ and $\S(0)=I_n$, we use (\ref{v1.11.4}),  (\ref{v1.11.7}), and the third
equality in  (\ref{v1.11}) to rewrite    (\ref{v1.18})
in the form
\begin{equation} \label{v1.24}
 \om_1( z)  =I_p-i\Lam_2(0)^*(\a - z I_n)^{-1}\Psi_2, \quad   \om_2( z)  =-\frac{1}{2}D-i\Lam_1(0)^*(\a - z I_n)^{-1}\Psi_2.
\end{equation}
The following procedure is a standard one in system theory.
Rewrite the first equality in (\ref{v1.17}) as
\begin{equation} \label{v1.24'}
i \Psi_2  \Lam_2(0)^*=\a -\g=(\a - z I_p)-(\g -z I_p)
\end{equation}
to check that
\[
(I_p-i\Lam_2(0)^*(\a - z I_n)^{-1}\Psi_2)(I_p+i\Lam_2(0)^*(\g - z I_n)^{-1}\Psi_2)=I_p.
\]
In other words, we have
\begin{equation} \label{v1.25}
 \om_1( z)  ^{-1}=I_p+i\Lam_2(0)^*(\g - z I_n)^{-1}\Psi_2.
\end{equation}
From (\ref{v1.24})-(\ref{v1.25}) it follows that
\begin{equation} \label{v1.26'}
i  \om_2( z)   \om_1( z)  ^{-1}=-\frac{i}{2}D+\Psi_1(0)^*(\g - z I_p)^{-1}\Psi_2,
\end{equation}
that is, the function given by (\ref{v1.16}) is a Weyl function.

Now, consider the matrix function $\wh \vp$, given by \eqref{c1.16}. If $D<0$,
we have $\vp=\wh \vp$ and  the case of $\vp$ is as above.
If the inequality $D<0$ does not hold,
we can introduce two diagonal matrices ${ P}_j={ P}_j^*={ P}_j^2$ 
($j=1,2$) such that
\begin{equation} \label{r4}
{P}_1+ { P}_2=I_p, \quad D \big({ P}_1- { P}_2\big)=|D|.
\end{equation}
Here ${P}_1$ has nonzero entry (entry equal to $1$), when
the  entry of $D$ in the same row and column is positive, and ${P}_2$ has nonzero entry (entry equal to $1$), when
the corresponding entry of $D$ is negative.
Let us show that for $\Im  z  \geq \|\a\|$ we have
\begin{equation} \label{r5}
\int_0^\infty[{P}_1 \quad {P}_2]
 \Om( z)  ^*w(x, z)  ^*H(x)w(x, z)   \Om( z)  \left[\begin{array}{c}
{P}_1 \\  {P}_2
\end{array}
\right] dx <\infty.
 \end{equation}
 It suffices to show that 
 \begin{equation} \label{r6}
\int_0^\infty[0 \quad {P}_2]
 \Om( z)  ^*w(x, z)  ^*H(x)w(x, z)   \Om( z)  \left[\begin{array}{c}
0 \\  {P}_2
\end{array}
\right] dx <\infty,
 \end{equation}
 which is immediately evident from (\ref{v1.22}), and that  
 \begin{equation} \label{r7}
\int_0^\infty[{P}_1 \quad 0]
 \Om( z)  ^*w(x, z)  ^*H(x)w(x, z)   \Om( z)  \left[\begin{array}{c}
{P}_1 \\  0
\end{array}
\right] dx <\infty.
 \end{equation}
 To prove (\ref{r7}) note that similar to (\ref{v1.21}) one obtains  
  \begin{equation} \label{r8}
[{P}_1 \quad 0]
 \Om( z)  ^*w(x, z)  ^*H(x)w(x, z)   \Om( z)  \left[\begin{array}{c}
{P}_1 \\  0
\end{array}
\right] =Q_1(x)^*Q_1(x), 
 \end{equation}
  \begin{equation} \label{r9}
 Q_1(x)=\big(D -i\Psi_1(x)^*\S(x)^{-1}(\a -  z  I_n)^{-1}\Psi_1(x)\big)e^{i  z  x D}{P}_1.
   \end{equation}
 It is an immediate consequence of the above that the entries of $De^{i  z  x D}{P}_1$ belong $L^2(0, \infty)$.
It follows from  (\ref{v1.13}) that the entries of $\Psi_1(x)e^{i  z  x D}{P}_1$ are bounded
for $\Im  z  \geq \|\a\|$.  Because of \eqref{c4} 
 the entries of $\Psi_1(x)^*\S(x)^{-1}$ belong $L^2(0, \infty)$.
Thus, in view of  (\ref{r9})  the entries of $Q_1$ belong $L^2(0, \infty)$ ($\Im  z  \geq \|\a\|$),
and using (\ref{r8}) we see that (\ref{r7}) is true. Hence (\ref{r5}) is also valid.
Substitute $\wh \vp$ of the form
 \begin{equation} \label{r10}
 \wh \vp( z)  =i \wh \om_2( z)  \wh \om_1( z)  ^{-1}, \quad
\left[\begin{array}{c}
\wh \om_1( z)   \\  \wh \om_2( z)  
\end{array}
\right] :=
  \Om( z)  \left[\begin{array}{c}
{P}_1 \\  {P}_2
\end{array}
\right] 
    \end{equation}
into   the left-hand side of (\ref{v1.6}).
    By (\ref{r5})  the function $\wh \vp$ of the form \eqref{r10}
    satisfies  (\ref{v1.6}) in the domain $\Im z >\|\a\|$.
    
According to  (\ref{v1.18}), (\ref{r4}), and (\ref{r10}) we have
 \begin{equation} \label{r11}
\left[\begin{array}{c}
\wh \om_1( z)   \\ \wh \om_2( z)  
\end{array}
\right] =w_{\a}(0,  z)  
\left[\begin{array}{c}
I_p \\  |D|/2
\end{array}
\right] .
    \end{equation} 
In view of    (\ref{r11}),  the formula  
\begin{align}&\label{c5}
i \wh \om_2( z)  \wh \om_1( z)  ^{-1}=\frac{i}{2}|D|+\wh \Psi_1(0)^*(\wh\g - z I_n)^{-1}\wh \Psi_2
\end{align}
 can be proven quite similar
to  (\ref{v1.26'}).    Equation \eqref{c1.16} is an immediate consequence of \eqref{r10} and \eqref{c5},
and so $\wh \vp$ given by \eqref{c1.16} satisfies \eqref{v1.6}.

To prove the uniqueness of the Weyl function for $D<0$, notice that 
\begin{eqnarray} \nonumber && 
H_0+cJ=\left[ \begin{array}{c} D/2 \\
2D^{-1}\big(\frac{1}{2}D+cI_p\big)
\end{array} \right]\left[ \begin{array}{lr} D/2  &
2\big(\frac{1}{2}D+cI_p\big)D^{-1}
\end{array} \right]
\\ \label{v1.26} &&
-\left[ \begin{array}{lr} 0 & 0 \\
0 & 4cD^{-1}\big(D+cI_p\big)D^{-1}
\end{array} \right].
\end{eqnarray}
When $D<0$ and the scalar $c>0$ is sufficiently small, 
formula (\ref{v1.26}) implies $H_0 \geq -cJ$. Hence, because of (\ref{v1.11.1})-(\ref{v1.11.3}) we have
\begin{equation} \label{v1.27}
H(x)\geq -cv_0(x)^*Jv_0(x)=-cJ.
\end{equation}
Using  (\ref{v1.3})  we obtain
\begin{equation} \label{v1.28}
-c \frac{d}{dx} \big(w(x, z)  ^*Jw(x, z)\big)  =ic(\ov z -  z)  w(x, z)  ^*H(x)w(x, z)   \geq 0, \quad  z  \in \BC_+.
\end{equation}
From (\ref{v1.27})  and (\ref{v1.28}) it follows that
\begin{equation} \label{v1.29}
w(x, z)  ^*H(x)w(x, z)   \geq -c w(x, z)  ^*Jw(x, z)  \geq -c w(0, z)  ^*Jw(0, z)  =-cJ.
\end{equation}
Formula (\ref{v1.29}) implies that
\begin{equation} \label{v1.30}
\int_0^r 
\left[ \begin{array}{lr} I_p& -I_p 
\end{array} \right]
w(x, z)  ^*H(x)w(x, z)  \left[ \begin{array}{c} I_p \\
-I_p
\end{array} \right] dx \geq 2crI_p.
\end{equation}
Denote by $L\subset \BC^m$ ($m=2p$) the subspace  of vectors $f$ such that
\[
\int_0^r 
f^*
w(x, z)  ^*H(x)w(x, z)  f dx<\infty.
\]
Inequality (\ref{v1.30}) implies $\dim L \leq p$,
which implies the uniqueness of the Weyl function.

Finally,   formulas  (\ref{v1.9}) and (\ref{v1.17})-(\ref{c1.17'}) imply   (\ref{v1.17''}).
\end{proof}
\begin{Rk} Representations like \eqref{v1.16} and \eqref{v1.17} are called
{\it realizations} in Control Theory.
\end{Rk}

\begin{Rk}   It follows from the proof of Proposition \ref{PnWe} that $\vp$ satisfies \eqref{v1.6}
in $\BC_+\setminus \big(\s(\a)\cup \s(\g)\big)$
and $\wh \vp$
satisfies \eqref{v1.6}  in  $\{z:\, \Im z> \|\a\|\} $.

\end{Rk}

\subsection{Inverse problem: explicit solution}\label{susinv}
The following theorem allows us to recover $H$ explicitly from the rational Weyl functions.
\begin{Tm}\label{Invpr} 
Let the diagonal $p \times p$ matrix $D$ be negative (i.e., let $D<0$). 
If  $\vp$ is a rational matrix function such that
\begin{equation} \label{r12}
\lim_{z \to \infty} \vp(z)=\frac{i}{2}|D|, \quad \Im \vp( z) \geq 0 \quad (z \in \BC_+),
 \end{equation}
then $\vp$ admits a realization 
\begin{equation} \label{r13}
\vp(z)=\frac{i}{2}|D|+\Psi_1(0)^*(\g - z I_n)^{-1}\Psi_2,  
\end{equation}
where $\Psi_1(0)$ and $\Psi_2$ are $n \times p$ matricies, $n\in \BN$,
and $n \times n$ matrix $\g$ satisfies the matrix identity
\begin{equation} \label{r14}
 \g-\g^* =i\big(\Psi_1(0) -\Psi_2\big) D^{-1}\big(\Psi_1(0) -\Psi_2\big) ^*.
\end{equation}
Moreover $\vp$ is the Weyl function of  the canonical system
determined (in the sense of Definition \ref{Dn1}) by the parameter matrices
\begin{align}\label{r15}& \a=\g+i\Psi_2\Lam_2(0)^*, \quad
\Lam_1(0)=\frac{1}{2}\big(\Psi_1(0)+\Psi_2\big), \quad
\Lam_2(0)=\big(\Psi_1(0)-\Psi_2\big)D^{-1}.
\end{align}
\end{Tm}
\begin{proof}.
Realization \eqref{r13} follows directly from Proposition 4.1 in \cite{ASAK} (see also
Theorem 5.2 in \cite{GKS6}). The identity \eqref{v1.9} follows from \eqref{r14}
and \eqref{r15}, that is, the requirement \eqref{v1.9} for parameter matrices
is fulfilled.
Moreover, \eqref{r15} implies relations \eqref{v1.17} and \eqref{v1.17'}.
Now, compare \eqref{v1.16} and \eqref{r13} to see that $\vp$
is the Weyl function by Proposition \ref{PnWe}.
\end{proof}

\section{Canonical system. General formulas} \label{CanG}
\setcounter{equation}{0}
In the Krein paper \cite{Kr} (see also \cite{AGKLS0, AGKLS})
Krein system was treated as a system generated by an accelerant
$k(x)=k(-x)^*$, where $k$ was  the kernel of operator $S$ (with difference kernel) of the form
\eqref{1.6}. In this section we assume that the diagonal matrix $D$ is negative (i.e., $D<0$).
In a similar way to \cite{Kr}, operators with $|D|$-difference kernel generate
other subclasses of canonical systems (see  \cite{ASAK, SaL2c, SaL4}
and references therein):
\begin{Pn}  \label{PnG}
Let $k(x)=\{ k_{ij}(x) \}_{i,j=1}^p$  be  
a $p \times p$ matrix function such that
\begin{align}\label{g0}&
k(x) \in L^2_{p\times p}(0,l), \quad k(-x)=k(x)^*,   \quad 
S_l, \, S_l^{-1}\in\{L^2_p(0,\,l), \, L^2_p(0,\,l)\}, \\
&
\label{g1}
S_l:=I+\int_0^l\{k_{ij}(d_jt-d_ix)\}_{i,j=1}^p\, \cdot \, dt>0, \quad \mathrm{for \, all}
\quad 0<l<\infty,
\end{align}
that is, operators $S_l$ determined via $k$ are positive,
bounded, and boundedly invertible.
Then operators $S_l$ admit a triangular factorization
\begin{align}&
\label{g2}
 S_l^{-1}=E_l^*E_l, \quad E_l=I+\int_0^xE(x,t)\, \cdot \, dt\in \{L^2_p(0,\,l), \, L^2_p(0,\,l)\},
\end{align}
where $E(x,t)$ is a Hilbert-Schmidt kernel. 

Moreover,
$k$ generates, in terms of $S$ and $E$, the canonical system  \eqref{v1.3}
such that  \eqref{c1} holds, and the Hamiltonian $H=\b^*\b$ of
the canonical system is given via
\begin{align}&
\label{g3}
\b(x)=\big(E_l\Pi\big)(x), \quad \Pi(x):=\Big[D\{s_{ij}(|d_i|x)\}_{i,j=1}^p \qquad I_p\Big], \\
& \label{g3'}
\quad s(x):=\frac{1}{2}I_p+|D|^{-1}\int_0^xk(t)dt.
\end{align}
This $H(x)$ is summable on all  finite intervals $(0, \, l)$.

The fundamental solution of the canonical system is given by
\begin{align}
 & \label{g4}
w(l,z)=I_{2p}+i z J\Pi_l^*S_l^{-1}(I -z A)^{-1}\Pi_l, \quad A,\Pi_l \in \{\BC^{2p}, \, L^2_p(0,\,l)\};
\\ &  \label{g5}A=A_l=iD\int_0^x \, \cdot \, dt; \quad \Pi_l g=\Pi(x)g, \quad g \in \BC^{2p},
\end{align}
where $\BC^{2p}$ is the $2p$-dimensional vector space.
\end{Pn}
\begin{proof}. According to \cite{SaLk1} (Ch. 6) the operators
$A, \, S_l, \, \Pi_l$ given by \eqref{g1} and \eqref{g5} satisfy the operator identity
\begin{align}\label{g6}
 &AS_l-S_lA^*=i \Pi_l J\Pi_l^*,
\end{align}
that is, they form an $S$-node \cite{SaL2, SaLk1, SaL2c, SaL4}. By Theorem 4.2.1
\cite{SaL4} there is a system corresponding to the family $A, \, S_l, \, \Pi_l$
($0<l<\infty$) of $S$-nodes, and $w$ given by \eqref{g4} is its fundamental solution.
Moreover, if $B(l):=\Pi_l^*S_l^{-1}\Pi_l$ is differentiable, this system is canonical, and
its Hamiltonian is given by
\begin{align}\label{g7}
H(l)=\frac{d}{dl}B(l)=\frac{d}{dl}\Pi_l^*S_l^{-1}\Pi_l.
\end{align}
The factorisation $S_l^{-1}=(I+E_+)(I+E_-)$, where $E_+$ ($E_-$)  is an upper (lower)
 triangular integral operator,  is clear from the factorization result on p. 184 in \cite{GoKrb}.
 Now, taking into account $S_l=S_l^*$  it is easy to derive \eqref{g2}
(see, for instance, Section 4 in \cite{ASAK}). According to \eqref{g2} and
 \eqref{g7} we get $H=\b^*\b$, where $\b$  is given by \eqref{g3}.
 Note that $E(x,t)$ does not depend on $l$ since the lower-upper
 factorisation of $S_l$ is unique, and the lower-upper
 factorisation of $S_l$ is unique because the lower-upper
 factorisation of $I$ is unique.
 
It remains to be shown that the equality
\eqref{c1} holds  for almost all values of $x$. For that purpose we multiply  \eqref{g6} by $E_l$
from the left  and by $E_l^*$
from the right. Then, taking into account \eqref{g2}  
and \eqref{g3} we get:
\begin{align}&
\label{g8}
E_lAE_l^{-1}-\big( E_lAE_l^{-1} \big)^*=i(E_l\Pi)J(E_l\Pi)^*=i\b(x)J\int_0^l\b(t)^*\,\cdot \, dt.
\end{align}
As $E_lAE_l^{-1}$ is a lower triangular operator,  formula \eqref{g8} implies
\begin{align}\label{g9}&
E_lAE_l^{-1}=i\b(x)J\int_0^x\b(t)^*\,\cdot \, dt.
 \end{align}
Denote the kernel of 
$E_l^{-1}$ by $\wh E$ and rewrite the equality \eqref{g9} in terms of the kernels
of integral operators:
\begin{align}\label{g10}
D&+\int_t^x\big(D\wh E(y,t)+E(x,y)D\big)dy+\int_t^xE(x,y)D\int_t^y\wh E(r,t)drdy
\\ \nn &
=\b(x)J\b(t)^* \quad (t \leq x).
 \end{align}
Equality \eqref{c1} follows from \eqref{g10}.
\end{proof}
\begin{Dn}\label{DnG} The class of  Hamiltonians $H(x)$
$(0<x<\infty)$ of canonical systems generated by
matrix functions $k$, which satisfy the conditions of  Proposition \ref{PnG},
is denoted by ${\cal H}(D)$. 
\end{Dn}
We next consider some fixed ${ H}\in {\cal H}(D)$ and turn our focus towards M\"obius (also called linear-fractional) transformations
\begin{align}       \label{g11}
\varphi (z, l)=i  \big({\cal W}_{11}( z){\cal P}_{1}( z)
+{\cal W}_{12}( z ){\cal P}_{2}( z )\big)\big({\cal W}_{21}( z ){\cal P}_{1}( z)
+{\cal W}_{22}( z ){\cal P}_{2}( z)\big)^{-1} ,
\end{align}
where
\begin{align}\label{g12}&
{\cal W}(l,z)=\{{\cal W}_{ij}(z)\}_{i,j=1}^2=w(l,\ov z)^*,
 \end{align}
and the pairs $\{\clp_1(z), \, \clp_2(z)\}$ are pairs of
$p \times p$ matrix functions, which are meromorphic in $\BC_+$ and satisfy
(excluding, possibly, a discrete set of points)
the following relations
\begin{align}\label{g13}&
\clp_1^*\clp_1+\clp_2^*\clp_2>0, \quad \clp_1^*\clp_2+\clp_2^*\clp_1 \geq 0,
\end{align}
that is, $\{\clp_1, \, \clp_2\}$ are pairs of nonsingular matrix functions with property-$J$.
The set of  M\"obius transformations \eqref{g11}, where the pairs $\{\clp_1, \, \clp_2\}$
vary and $l$ is fixed, is denoted by $\cln_l$.
It follows from Statement 3 in \cite{SaA0} and
interpolation results in \cite{SaL2} (see also their
formulation in Theorem 1 from  \cite{SaA0} or Chapter 1 in \cite{SaL30})
that  there is a unique matrix function
\begin{align}\label{g14}&
\displaystyle{\vp= \bigcap_{l>0}\cln_l},
\end{align}
and this $\vp$ admits  representation
\begin{align}\label{g15}&
\vp(z) =zD^2\int_0^{\infty}e^{-izxD}\{s_{ij}(|d_i|x)\}_{i,j=1}^pdx\quad (z \in \BC_+),
\end{align}
where $s$ is given by \eqref{g3'}. Moreover, according to Theorem 1 from  \cite{SaA0}
the scalar products $\langle S_l f,f\rangle \,$ ($0<l<\infty$) admit representation
\begin{align}\label{g15'}&
\langle S_l f,f\rangle =\int_{-\infty}^{\infty}\Big(\int_0^le^{-itxD}f(x)dx\Big)^*
d\tau(t)\Big(\int_0^le^{-ityD}f(y)dy\Big).
\end{align}
Here $\tau$ is the $p \times p$ non-decreasing matrix function from Herglotz representation of $\vp$:
\begin{equation}       \label{a1}
\varphi ( z )= \mu z + \nu +
\int_{- \infty }^{ \infty}
\Big( \frac{1}{t- z } - \frac{t}{1+t^{2}}\Big)d \tau (t)= \mu z + \nu +
\int_{- \infty }^{ \infty}
 \frac{1+tz}{t- z }  \frac{d \tau (t)}{1+t^2},
\end{equation}
where  $\mu \geq 0\, $ ($\mu=0$ in our special case), $\nu=\nu^*$,  and
\begin{equation}       \label{a2}
\int_{-\infty}^{\infty}(1+t^2)^{-1}d\tau(t)<\infty.
\end{equation}

Note that according to \eqref{v1.3} and \eqref{g12}  we have
\begin{align}\label{g16}&
\frac{d}{dx}\big(w(x,\ov z)^*Jw(x,z)\big)=0, \quad {\mathrm{and \, \, so}} \quad
{\cal W}(l, z)=Jw(l,z)^{-1}J.
\end{align}
From \eqref{g16}  we get
\begin{align}\label{g18}&
\big({\cal W}(l, z)^*\big)^{-1}J{\cal W}(l, z)^{-1}=Jw(l, z)^*Jw(l,z)J.
\end{align}
Moreover, taking into account \eqref{g4} and \eqref{g6}
we obtain
\begin{align} \label{g20}&
w(l, z)^*Jw(l,z)=J+i(z - \ov z)\Pi_l^*(I-\ov z A^*)^{-1}S_l^{-1}(I-zA)^{-1}\Pi_l.
\end{align}
Using \eqref{g3} and  \eqref{g5}, separate $\Pi_l$ into two blocks $\Phi_k \in \{\BC^p, \, L^2_p(0, \, l)\}$:
\begin{align} \label{g21}&
\Phi_1h=D\{s_{ij}(|d_i|x)\}_{i,j=1}^ph,  \quad \Phi_2h\equiv h \quad (h \in \BC^p);
\quad \Pi_l=[\Phi_1 \quad \Phi_2].
\end{align}
We omit index "$l$" in $A$ and $\Phi_k$ as the choice of   "$l$" is clear from the context.
Finally, formulas \eqref{g18}-\eqref{g21} imply
\begin{align} \label{g22}r(l,z):&=-[I_p \quad 0]\big({\cal W}(l, z)^*\big)^{-1}J{\cal W}(l, z)^{-1}
\left[
\begin{array}{c} I_p 
 \\ 0
\end{array}
\right]
\\ \nn &=i(\ov z - z)\Phi_2^*(I-\ov z A^*)^{-1}S_l^{-1}(I-zA)^{-1}\Phi_2>0 \quad (z\in \BC_+).
\end{align}
By  \eqref{g14} and \eqref{g22} the conditions of Proposition 9.1.5 from \cite{SaL4} are fulfilled,
and hence $\vp$ satisfies \eqref{v1.6}. 
\begin{Tm} \label{DirG} Let $H\in {\cal H}(D)$ $\,(D<0)$. Then the unique
Weyl function of the corresponding canonical system \eqref{v1.3}
admits the representation 
\begin{align}\label{g22'}&
\vp(z) =-zD\int_0^{\infty}e^{izx}s(x)dx \qquad (z \in \BC_+),
\end{align}
where $s$ is given by  \eqref{g3'}.
\end{Tm}
\begin{proof}. 
Representation \eqref{g22'} is clear from \eqref{g15}.

It remains only to prove that $\vp$ satisfying \eqref{v1.6} is unique.
Here, it would be of interest to give a proof connected with $S$-nodes to compare
with the proof of uniqueness from Proposition \ref{PnWe}.
For this purpose let us consider $\Phi_2^*(I-\ov z A^*)^{-1}S_l^{-1}(I-zA)^{-1}\Phi_2$.

First, fix some $a\in \BR_+$, put $l=na\,$ ($n\in \BN$),  and consider the scalar product   in \eqref{g15'}.
It is clear that
\begin{align} \label{g23}&
\langle S_{na} f,f\rangle \leq p \sum_{j=1}^p \langle S_{na} f_j\E_j,f_j\E_j\rangle, \quad
f=:\{f_j\}_{j=1}^p, \quad \E_j:=\{\delta_{ij}\}_{i=1}^p.
\end{align}
Formula \eqref{g15'} and the equality below
\begin{align} \label{g24}&
\int_0^{na}e^{-ityD}\big(f_j(y)\E_j\big)dy=\sum_{k=1}^n e^{-i(k-1)t ad_j}\int_0^ae^{-ityd_j}f_{jk}(y)dy\, \E_j,
\\ \nn &
f_{jk}(y):=f_j\big(y+(k-1)a\big) \in L^2(0,a),
\end{align}
imply that
\begin{align} \label{g25}&
\langle S_{na} f_j\E_j,f_j\E_j\rangle \leq n \sum_{k=1}^n
\int_{-\infty}^{\infty}\Big|\int_0^ae^{-ityd_j}f_{jk}(y)dy\Big|^2\E_j^*d\tau(t)\E_j.
\end{align}
It follows from  \eqref{g15'}   and  \eqref{g25}  that
\begin{align} \label{g26}&
\langle S_{na} f_j\E_j,f_j\E_j\rangle \leq n \|S_a\|\sum_{k=1}^n\|f_{jk}\|^2=n \|S_a\|\|f_{j}\|^2.
\end{align}
By \eqref{g23} and \eqref{g26} we have $\langle S_{na} f,f\rangle \leq pn \|S_a\|\|f\|^2$,
 and so we get
\begin{align} \label{g27}&
\|S_{na}\| \leq pn \|S_a\|, \quad {\mathrm{i.e.,}}\quad S_{na}^{-1}\geq \big(pn \| S_a\|\big)^{-1}I.
\end{align}
One can easily check that $(I-zA)^{-1}\Phi_2=e^{izxD}$. Thus, using \eqref{g27} we
derive
\begin{align} \label{g28}&
\Phi_2^*(I-\ov z A^*)^{-1}S_{na}^{-1}(I-zA)^{-1}\Phi_2 \geq c(n)I_p, \quad c(n) \to +\infty \, (n\to \infty).
\end{align}

Now, let $\vp$ and $\wt \vp$, such that $\psi(z):=\vp(z)- \wt \vp(z)\not=0$ for some $z\in \BC_+$,
satisfy \eqref{v1.6}.  This implies that
\begin{align} \label{g29}&
\int_0^\infty
\left[ \begin{array}{lr}
0 &i \psi(z)^* \end{array} \right]
  w(x,  z)  ^*H(x)w(x,  z)  
\left[ \begin{array}{c}
0 \\ -i\psi(z)  \end{array} \right]
dx  <\infty.
\end{align}
According to \eqref{v1.3} and similar to \eqref{v1.28} we have 
\[
\frac{d}{dx}\Big(w(x,z)^*Jw(x,z)\Big)=i(z-\ov z)w(x,z)^*H(x)w(x,z).
\]
Hence, we rewrite \eqref{g29} as
\begin{align} \label{g30}&
\sup_{l<\infty}\left(i(\ov z -z)^{-1}\left[ \begin{array}{lr}
0 &i \psi(z)^* \end{array} \right]
  w(l,  z)  ^*Jw(l,  z)  
\left[ \begin{array}{c}
0 \\ -i\psi(z)  \end{array} \right]\right)
 <\infty .
\end{align}
In view of \eqref{g20} inequality \eqref{g30} implies
\begin{align} \nn&
\sup_{n\in \BN}\Big(\psi(z)^* 
 \Phi_2^*(I-\ov z A^*)^{-1}S_{na}^{-1}(I-zA)^{-1}\Phi_2
\psi(z)  \Big)
 <\infty ,
\end{align}
which contradicts \eqref{g28}.
\end{proof}
Theorem \ref{DirG} yields a solution of the inverse problem directly via Weyl function.
\begin{Tm} \label{InvG} Let $\vp$ be the Weyl function of the
canonical system such that its Hamiltonian $H\in {\cal H}(D)$ $\,(D<0)$. Then we
have
\begin{align}&\label{g31}
e^{-\eta x}s(x) \in \Big(L^1_{p\times p}(0, \, \infty) \bigcap L^2_{p\times p}(0, \, \infty)\Big)
\end{align}
for any $\eta >0$, and the matrix function  $k(x)$ $\,(x>0)$ is recovered via Fourier transform
\begin{align}&\label{g32}
k(x)= \frac{1}{2\pi} \frac{d}{dx}\Big(e^{\eta x}
{\mathrm{l.i.m.}}_{a\to \infty} \int_{- a}^{a}e^{-i \zeta x}\frac{\vp(\zeta +  i \eta)}{\zeta +  i \eta}
d\zeta \Big).
\end{align}
The Hamiltonian $H$ is recovered from $k$ as in Proposition \ref{PnG} or, equivalently,
using \eqref{g7}.
\end{Tm}
\begin{proof}. To prove \eqref{g31} recall the first inequality in  \eqref{g27},
which holds for all $a \in \BR_+$. It follows that
\begin{align}&\label{g33}
\|S_l\| \leq  \|S_{\lceil l \rceil}\| \leq p(l+1)\|S_1\| \leq C_1(l+1), \quad C_1 \in \BR_+.
\end{align}
We next represent the operator $\Phi_1=\Phi_{1,l}$ given by \eqref{g21} in the form
\begin{align}&\label{g34}
\Phi_1=\Up_l(z)+i\Phi_2\vp(z)=(I-zA)(I-zA)^{-1}\Up_l(z)+i\Phi_2\vp(z), 
\\ &
\nn \Up_l(z):=\Phi_1-i\Phi_2\vp(z).
\end{align}
Now, we need  inequality (22) from \cite{SaA0}:
\begin{align}&\label{g35}
\| (I-zA)^{-1}\Up_l(z)\|\leq \|S_l\|^{1/2}\sqrt{\big\|\big(\vp(z)-\vp(z)^*\big)/(z- \ov z)\big\|}, \quad
z \in \BC_+.
\end{align}
Note that we can fix in \eqref{g34} and \eqref{g35} any $z\in \BC_+$.
Taking into account relations $\|\Phi_2\|= \sqrt{l}$ and $\|A\|<C_2l$,
we derive from \eqref{g33}-\eqref{g35} that
\begin{align}&\label{g36}
\|\Phi_1\| \leq C_3(l+1)^{3/2}.
\end{align}
Relation \eqref{g31} follows from \eqref{g21} and \eqref{g36}.

Rewrite \eqref{g22'} as
\begin{align}\label{g37}&
|D|^{-1}\vp(\zeta+i\eta)/(\zeta+i\eta) =\int_0^{\infty}e^{i\zeta x}e^{-\eta x}s(x)dx \qquad (\zeta\in \BR, \, \, \eta>0).
\end{align}
In view of the Plancherel theorem formulas \eqref{g31} and  \eqref{g37} imply
\begin{align}&\label{g38}
|D|^{-1}\vp(\zeta+i\eta)/(\zeta+i\eta) ={\mathrm{l.i.m.}}_{a\to \infty}\int_0^{a}e^{i\zeta x}e^{-\eta x}s(x)dx, \\
&\label{g39}
e^{-\eta x}s(x)=\frac{1}{2\pi}{\mathrm{l.i.m.}}_{a\to \infty}\int_{-a}^{a}e^{-i\zeta x}|D|^{-1}\vp(\zeta+i\eta)/(\zeta+i\eta)d\zeta 
\end{align}
for all fixed $\eta >0$. Note that \eqref{g39} holds for $x>0$, that is, l.i.m. in \eqref{g39}
is considered in $L^2_{p\times p}(0, \,\infty)$. Finally, equalities \eqref{g3'} and  \eqref{g39} yield \eqref{g32}.
\end{proof}
Using Theorems \ref{DirG} and \ref{InvG} we obtain our next Borg-Marchenko-type result.
\begin{Tm}  \label{BM} Let $\vp$  and $\wt \vp$  be the Weyl functions of the
canonical systems with Hamiltonians $H$ and $\wt H$, respectively, where $H, \, \wt H \, \in \, {\cal H}(D)$ $\,(D<0)$. 
Let the equalities
\begin{align}&\label{g40}
e^{-iza}\big(\vp(z)-\wt \vp(z)\big)=O(1) \quad (z\to \infty)
\end{align}
hold for all $a<l$ on some ray $\Im z/\Re z=c$ in $\BC_+$. Then we have
\begin{align}&\label{g41}
H(x) \equiv \wt H(x) \quad {\mathrm{for}} \quad 0<x<l/d, \quad d:= \max_{1\leq k \leq p}|d_k|.
\end{align}
\end{Tm}
\begin{proof}.
We mark functions corresponding  to $\wt \vp$ with a tilde, e.g., $\wt H, \, \wt s$.
Consider the entire matrix function
 \begin{align}&\label{g42}
\om(z):=\int_0^a      e^{iz(x-a)}            \big(s(x)- \wt s(x)\big)dx.
\end{align}
It is clear that $\| \om(z)\|$ is bounded in the closed lower semi-plane and tends to zero on some ray
$\Im z/\Re z=c_1$ there. To estimate
$\om$ on the ray $\Im z/\Re z=c$ in $\BC_+$ we use \eqref{g22'} and write
 \begin{align}&\label{g43}
\om(z)=e^{-iza}(z|D|)^{-1}\big(\vp(z)-\wt \vp(z)\big)-\int_a^{\infty} e^{iz(x-a)}            \big(s(x)- \wt s(x)\big)dx.
\end{align}
From \eqref{g31},  \eqref{g40}, and \eqref{g43} one gets the boundedness of $\| \om(z)\|$ on the ray $\Im z/\Re z=c$.
Hence, according to the Phragmen-Lindel\"of theorem $\| \om(z)\|$ is bounded in $\BC_+$,
and thus also in $\BC$. As $\| \om(z)\|$ tends also to zero on some ray, we see that $\om(z) \equiv 0$, that is,
$s(x) \equiv \wt s(x)$ for $x <a$. Since this fact is true for all $a <l$ the equalities
 \begin{align}&\label{g44}
s(x) \equiv \wt s(x), \quad k(x) \equiv \wt k(x) \quad (0<x <l)
\end{align}
follow. Now, take into account \eqref{g1},  \eqref{g3}, \eqref{g7},  and \eqref{g44} to obtain
\eqref{g41}.  
\end{proof}

\section{Interpolation of the Weyl function} \label{Intrp}
As can also be seen from the cases treated in Sections \ref{saDir}-\ref{CanG}
Weyl functions contain important information about systems, and their interpolation
is of particular interest. One of the possible approaches is interpolation by rational Weyl functions.
Recall that an inverse problem for rational Weyl functions of a subclass of canonical systems
was solved explicitly in Subsection \ref{susinv}. In this section we shall consider another
approach, namely an approach  in the spirit of \cite{RT}.
Theorem 4 from \cite{RT} (after reformulation for the matrix case and for the upper semi-plane) has the following form.
\begin{Tm}\label{Tuan} Let matrix function $F$  admit the representation
\begin{equation} \label{1.14}
F(z)= \int_0^{\infty}e^{izx}f(x)dx, \quad e^{-\eta x}f(x)\in L^1_{p \times p}(0, \, \infty) \,\, {\mathrm{for \,\, all}}\,\, \eta>0.
\end{equation} 
Then for any $\ve>0$ and $\Im z>\frac{1}{2}+\ve$ we have
\begin{align} \label{1.15}&
F(z)=\sum_{n=0}^{\infty}c_n\big(z+\frac{i}{2}-i\ve\big)\sum_{q=0}^n a_{nq}F(iq+i\ve),
 \end{align} 
where
\begin{align} \label{1.16}&
a_{nq}=\frac{(-1)^q(n+q)!}{(q!)^2(n-q)!}, \quad c_n(\la)=\frac{(2n+1
)\prod_{q=1}^n(q-\frac{1}{2}+i\la)}
{\prod_{q=0}^n(q+\frac{1}{2}-i\la)}.
 \end{align} 
The series in \eqref{1.15} converges uniformly  on any compact subset in the semi-plane $\Im z>\frac{1}{2}+\ve$ and we get
\begin{equation} \label{1.17}
\| F(z)-\sum_{n=0}^{N}c_n\big(z+\frac{i}{2}-i\ve\big)\sum_{q=0}^n a_{nq}F(iq+i\ve)\|=
O\Big(N^{\frac{1}{2}+\ve-\Im z}\Big), \quad N\to \infty.
 \end{equation}  
\end{Tm}
Our next theorem is an analogue for the Dirac system case of the main Theorem 5
from \cite{RT}.
\begin{Tm}\label{IntTm} Let $\vp$ be the Weyl function of system \eqref{1.1}, where $v$ satisfies
\eqref{1.10'}. Then for any $\ve>0$ and $\Im z>\frac{1}{2}+\ve$ the matrix function $\vp$ admits the representation
 \begin{align} \label{1.18}&
\vp(z)=-z^2\sum_{n=0}^{\infty}c_n\big(z+\frac{i}{2}-i\ve\big)\sum_{q=0}^n (q+\ve)^{-2}a_{nq}\vp(iq+i\ve).
 \end{align} 
Moreover,  for  $\Im z>\frac{1}{2}+\ve$ we get
 \begin{align}\nonumber &
\|\vp(z)+z^2\sum_{n=0}^{N}c_n\big(z+\frac{i}{2}-i\ve\big)\sum_{q=0}^n (q+\ve)^{-2} a_{nq}\vp(iq+i\ve)\|
\\ &  \label{1.19}
=
O\Big(N^{\frac{1}{2}+\ve-\Im z}\Big), \quad N\to \infty.
 \end{align} 
\end{Tm}
\begin{proof}. The theorem's statement follows from formula \eqref{1.12} and  Theorem \ref{Tuan}.
 \end{proof}
\begin{Rk}
 One can use various transformations of $\vp$ to change the interpolation set $\{z: \, z= iq+i\ve, \,\, q\in \BN_0\}$.
The simplest example is $\vp(z+z_0)$, where $z_0 \in (\BC_+\cup \BR)$. Noticing that  $\vp(z)/z^2$ satisfies the conditions
of  Theorem \ref{Tuan}, we see that $\vp(z+z_0)/(z+z_0)^2$ (as a function of $z$) also satisfies  these conditions. 
Hence, we obtain
 \begin{align}\nonumber
\vp(z+z_0)=&-(z+z_0)^2\sum_{n=0}^{\infty}c_ n\big(z+\frac{i}{2}-i\ve\big)
\\  \label{1.20}& \times
\sum_{q=0}^n (q+\ve-iz_0)^{-2}a_{nq}\vp(z_0+iq+i\ve).
 \end{align} 
\end{Rk}

Clearly, representation \eqref{g22'} enables us to apply Theorem \ref{Tuan} to
Weyl functions of canonical systems.

In fact,  interpolation formulas  of  \eqref{1.18}  type are true for the Weyl functions of the  much
wider class of Dirac systems as well as  for other classical and non-classical systems
because the corresponding Weyl functions are either bounded or belong to the
Nevanlinna  (also called Herglotz)  class, that is, $\Im \vp(z) \geq 0$.

Indeed,
according to Proposition 4.2 in \cite{SaA2} there is a unique Weyl function $\vp$
associated with a Dirac system on $[0, \infty)$ with  a locally summable potential $v$. Moreover,
this $\vp$ belongs to the Nevanlinna   class
and therefore admits a Herglotz representation \eqref{a1}.

\begin{Pn}\label{Herg} Suppose $\vp(z)$ $(z \in \BC_+)$ belongs to the Nevanlinna class.
Then the matrix functions $(z+i \delta)^{-2} \vp(z+i \delta)$ $(\delta>0)$ satisfy the
conditions of  Theorem \ref{Tuan}.
\end{Pn}
\begin{proof}.
Put $z+i \delta = \xi + i \eta$ ($\xi \in \BR$, $\, \eta >\delta$).
Then, in view of \eqref{a1} and \eqref{a2} (and for some $C_1, \, C_2>0$) we get the inequality
\begin{align}&\nn
\left(\int_{-\infty}^{\infty}|\xi + i \eta|^{-4}|h^*\vp(\xi + i \eta)h|^2  d \xi \right)^{1/2}
\\ \nn & 
\leq
\left(\int_{-\infty}^{\infty}\frac{1}{\xi^2 +  \eta^{2}}
\left|\int_{-\infty}^{\infty}\frac{\xi d\tau_h(t)}{(t-\xi-i\eta)(1+t^2)}
\right|^2d\xi \right)^{1/2}+C_1
\\ \nn & \leq \left(\int_{-\infty}^{\infty} 
\int_{-\infty}^{\infty}
\int_{-\infty}^{\infty} \Big((t-\xi)^2+\eta^2\Big)^{-\frac{1}{2}}\Big((\zeta-\xi)^2+\eta^2\Big)^{-\frac{1}{2}}d\xi
\frac{ d\tau_h(\zeta)}{1+\zeta^2}\frac{ d\tau_h(t)}{1+t^2}\right)^{\frac{1}{2}}
\\ \label{a3} &+C_1
\leq \int_{-\infty}^{\infty} \frac{ d\tau_h(t)}{1+t^2} \left(\int_{-\infty}^{\infty}
\frac{ d\xi}{\xi^2+\delta^2} \right)^{1/2} +C_1\leq C_2,
\end{align}
where $h\in \BC^p$, $h^*h \leq 1$, $\tau_h(t)=h^*\tau(t) h$, and
$C_2$ does not depend on $h$ and on $\eta>\delta$. Therefore,  the norms of entries
of $(\xi+i\eta)^{-2}\vp(\xi+i\eta)$ in $L^2(-\infty, \, \infty)$ are uniformly bounded for
$\eta>\delta $.
Thus, we can apply Theorem V from \cite{WP}
and derive
\begin{align}\label{a4}
 &(z+i \delta)^{-2} \vp(z+i \delta)={\mathrm{l.i.m.}}_{a \to \infty}\int_0^ae^{izx}f(x)dx \quad (\Im z >0), 
\end{align}
where $f \in L^2_{p \times p}(0,\, \infty)$. Representation \eqref{1.14} follows from \eqref{a4}
and from the fact that  $(z+i \delta)^{-2} \vp(z+i \delta)$ is analytic.
\end{proof}
\begin{Rk}\label{kappa}
The integral representation \eqref{1.14} is true also for  functions $\vp$ from the generalized
Nevanlinna class $N_{\vk}$ under additional conditions
\begin{align}& \nn
\lim_{\eta \to \infty}\vp(i \eta)=0, \quad \lim_{\eta \to \infty}\eta |\Im \vp(i \eta)|<\infty
\end{align}
in some semi-planes $\Im z >h_{\vp}$ \cite{KL}. Weyl functions from 
the generalized
Nevanlinna class appear, for instance, in \cite{Ka, KST, LW, RS} (see also references therein).
 \end{Rk}

{\bf Acknowledgement.}
The work of A.L. Sakhnovich was supported by the Austrian Science Fund (FWF) under
Grant  no. Y330.

\end{document}